\documentclass[11pt,lefteqn]{article}

\usepackage{amsmath,amsfonts,amsthm,amssymb,graphics,epsfig,color}
\newcommand{\R}{\mathbb R}
\newcommand{\C}{\mathbb C}
\newcommand{\N}{\mathbb N}

\renewcommand{\Re}{\mathop {\rm Re}\nolimits}
\renewcommand{\Im}{\mathop {\rm Im}\nolimits}

\newcommand*{\clos}[1]{\overline{#1}}

\newtheorem{theorem}{Theorem}
\newtheorem{lemma}[theorem]{Lemma}

\begin{document}
\date {\today}
\title{Update to:  {\em Spectral sets: numerical range and beyond}}

\author{Michel Crouzeix\footnote{Univ.\,Rennes, CNRS, IRMAR\,-\,UMR\,6625, F-35000 Rennes, France.
email: michel.crouzeix@univ-rennes1.fr},  Anne Greenbaum\footnote{University of Washington,
Applied Math Dept., Box 353925, Seattle, WA 98195.  email:  greenbau@uw.edu}
}

\maketitle

\begin{abstract} 
We use the results of Lorist and Schwenninger, {\em A solution to Crouzeix's
conjecture}, arXiv:2608.03841v2, https://arxiv.org/abs/2608.03841,
to update the 2019 paper of Crouzeix and Greenbaum, {\em Spectral
sets:  Numerical range and beyond}, SIAM J.~Matrix Anal.~Appl.,
40(3):1087-1101.  For all regions considered in the 2019 paper,
we are able to show that the region is a $K$-spectral set with
a bound on $K$ that is smaller than that established in the 2019 paper.
In particular, we show that various annular regions
are $2$-spectral sets and that a more general convex region with a circular hole
or cutout is a $4$-spectral set.
\end{abstract}

\paragraph{2000 Mathematical subject classifications\,:}47A25 ; 47A30

\noindent{\bf Keywords\,:}{ numerical range, spectral set}

\section{Introduction}
Let us consider a closed subset $X \subset\C$ of the complex plane and a bounded linear operator $A$ in a complex Hilbert space $(H, \langle, \rangle, \| \,\|)$. We will say that $X$ is a $K$-spectral set for 
$A$ if the spectrum of $A$ is contained in $X$ and if the following inequality
\begin{equation}
\label{eq1}
 \|f(A)\|\leq K\sup_{z\in X}|f(z)|,
\end{equation}
holds for all rational functions $f$ bounded in $X$. Note that $f(A)$ is naturally defined for such $f$ since, being bounded, $f$  has no pole in $X$. Let us denote by $\mathcal A(X)$ the set of
uniform limits in $X$ of  bounded rational functions; then,  by continuity, this inequality allows us to define $f(A)$ for $f\in \mathcal A (X)$ and inequality \eqref{eq1} still holds.

If $X$ is the closure of the numerical range $W(A)$, defined by
\begin{equation}
W(A) := \{ \langle Aq,q \rangle : q \in H ,~\| q \| = 1 \} , \label{numericalrange}
\end{equation}
then it is now known \cite{jin,schw} that $X$ is a $2$-spectral set for $A$.  In this paper, we extend this result
to show that other regions in the complex plane are $K$-spectral sets.  In particular, we show that various annular regions
are $2$-spectral sets and that a more general convex region with a circular hole
or cutout is a $4$-spectral set.  
We demonstrate how these results can be used to give bounds on the convergence rate of the GMRES algorithm for solving linear 
systems and on that of rational Krylov subspace methods for approximating $f(A)b$, where $A$ is a square matrix, $b$ is a given vector, 
and $f$ is a function that can be uniformly approximated on such a region by rational functions with poles outside the region.

Now we consider a bounded open subset $\Omega\subset \C$; we assume that its boundary $\partial \Omega$ is rectifiable and has a finite number of connected components;
  more precisely $\partial \Omega=\{\sigma (s)\,:s\in \partial \omega\}$ where $\partial \omega\subset\R$ is a finite union of disjoint segments $[a_j,b_j]$, $\sigma (a_j)=\sigma (b_j)$ and
  $s$ is the arc length of $\sigma (s)$ with a positive orientation.
 Then, if $A$ is a bounded linear operator with spectrum Sp$(A)$ contained in $\Omega$, it follows from the Cauchy formula that \eqref{eq1} holds with $X$ being the closure of $\Omega$ and $K=\frac{1}{2\pi }\int_{\partial \Omega}\|(\sigma I{-}A)^{-1}\|\,|d\sigma |$. 
But, this estimate is often very pessimistic, and we are looking for a better one. 
For that, we start with a rational function $f$ (bounded in $\Omega$) and we will consider the Cauchy formulae (for $z\in \Omega$)
\[
f(z)=\frac1{2\pi i}\int_{\partial \Omega}f(\sigma)\frac{d\sigma }{\sigma -z}, \quad f(A)=\frac1{2\pi i}\int_{\partial \Omega}f(\sigma)(\sigma I{-}A)^{-1}d\sigma .
\]
We will also consider the transforms of $f$ by the double layer potential kernel
\begin{equation}\label{eq3}
S(f,z):=\int_{\partial \omega}f(\sigma(s))\mu (\sigma(s),z)\,ds ,\quad S=S(f,A):=\int_{\partial \omega}f(\sigma (s))\mu (\sigma (s),A)\,ds.
\end{equation}
Here $\mu $ denotes the kernel\footnote{Note that $\mu $ is twice the usual kernel associated to the double layer potential.} given by 
\begin{align}\label{eq4}
\mu (\sigma(s),z):=\frac1\pi \frac{d\arg(\sigma (s){-}z)}{ds}=\frac{1}{2\pi i}\Big(\frac{\sigma '(s)}{\sigma(s) -z}-\frac{\clos{\sigma '(s)}}{\clos{\sigma(s)} -\bar z}\Big),\\
\label{eq5}M(s):=\mu (\sigma(s) ,A)=
\frac{1}{2\pi i}\big(\sigma '(s)(\sigma(s) I{-}A)^{-1}-\clos{\sigma '(s)}(\clos{\sigma(s)} I {-}A^*)^{-1}\big).
\end{align}
Note that $\sigma '(s)$ exists and $|\sigma '(s)|=1$ for almost every $s\in\partial \omega$ since $s$ is an arc length;
 these relations also assume $ \sigma (s)\neq z$, thus
 $\mu (\sigma(s),z)$ is defined a.e. $s\in \partial \omega$ and $\mu (\sigma(s),A)$ 
 is defined for almost all $s$ such that $\sigma (s)$ does not belong to the spectrum of $A$ . 
Note also that $\mu ( \sigma (s) , z)$ is real-valued and $\mu ( \sigma (s), A)$ is self-adjoint. 

Note also that, if we choose the constant function $f=1$, then $f(A)=I$,
\begin{equation}
\int_{\partial \omega}\mu (\sigma(s) ,z)\,ds=S(1,z)=2,\quad\text{if  }z\in\Omega\quad\text{and}\quad
\int_{\partial \omega}\mu (\sigma(s) ,A)\,ds=S(1,A)=2I. \label{S1zA}
\end{equation}

The following theorem is essentially the generalized version of Lemma 1 and Theorem 3 in
\cite[Remark 2 (ii)]{schw}.

\begin{theorem}\label{th}
Assume that  Sp$(A)\subset\Omega$ and there exists a bounded measurable real valued function $\lambda $ on $\partial \omega$ such that
\[
N(s):=\mu (\sigma (s),A){-}\lambda (s)\,I \quad\textrm{is positive semi-definite for almost all }s\in\partial \omega,
\]
then $\clos\Omega$ is a $K$-spectral set for the operator $A$ with a constant 
 \[
 K=\max\big(\rho ,1{+}\sqrt{\rho/2} \big)\quad \textrm{where  }\rho =\Big\|\int_{\partial \omega} N(s)\,ds\Big\|=2-\int_{\partial \omega}\lambda (s)\,ds.
 \]
\end{theorem}

\section{Proof of Theorem \ref{th}}

We get from Cauchy formula for a rational function $g$ bounded in $\Omega$  (writing $\sigma $ for $\sigma (s)$)
\begin{align*}
g(A)&=\frac{1}{2\pi i}\int_{\partial \Omega }g(\sigma )(\sigma I{-}A)^{-1}\,d\sigma =\int_{\partial \omega} g(\sigma ) M(s) \,ds +\frac1{2\pi i}\int_{\partial \Omega} g(\sigma )(\bar\sigma I{-}A^*)^{-1}\,d \bar{\sigma}\\
&=\int_{\partial \omega} g(\sigma ) N(s) \,ds-\int_{\partial \omega}  g(\sigma ) \lambda (s)I \,ds+\frac1{2\pi i}\int_{\partial \Omega} g(\sigma )(\bar\sigma I{-}A^*)^{-1}\,d \bar{\sigma}.
\end{align*}  
We will use the notation
\[
E(g):=\frac1{2\pi i}\int_{\partial \Omega} \overline{g(\sigma )}(\sigma I{-}A)^{-1}\,d\sigma
-\int_{\partial \omega}\overline{g(\sigma )}\lambda (s)I\,ds
=\int_{\partial \omega} \overline{g(\sigma )} N(s) \,ds - g(A)^*,
\]
the property $\int_{\partial \omega} N(s)\,ds =\rho \,I$, the semi definite positivity of $N(s)$
 and Cauchy-Schwarz inequalities
\begin{align*}
\int_{\partial \omega }|\langle N(s)u(\sigma ),v\rangle| \,ds&\leq \int_{\partial \omega}\langle N(s)u(\sigma ),u(\sigma )\rangle^{1/2} \langle N(s) v,v\rangle^{1/2} ds\\
&\leq \sqrt\rho \,\Big(\int_{\partial \omega}\langle N(s)u(\sigma ),u(\sigma )\rangle\,ds\Big)^{1/2}\,\|v\|.
\end{align*}

To prove the theorem, we argue ab absurdo. We consider rational functions $f$ bounded by 1 in $\Omega$, we set $T=f(A)$, $\kappa =\|T\|$ and we assume that $\kappa >\max \{ \rho , 1 + \sqrt{\rho /2} \}$; we choose a unit vector
$x$ such that $T^*Tx=\kappa ^2x$. (If $T$ is finite dimensional it is clear that such a vector $x$
exists; for the justification when $T$ is infinite dimensional, see \cite{schw}.)  We also set $u(\sigma )=\overline{f(\sigma )}Tx{-}\kappa x$ and
\begin{align*}
m_n&=\Re \langle E(f^n)T^nx,x\rangle=\Re\langle E(f^n)x,T^{*n}x\rangle\\
&=\Re \int_{\partial \omega } \overline{f(\sigma )}^n\langle N(s)x,T^{*n}x\rangle\,ds-\|T^{*n}x\|^2.
\end{align*}
We have used that $T$ (since this is the case for  $A$) commutes with $E(g)$. We have also
\begin{align*}
m_{n+1}&=\Re \langle E(f^{n+1})Tx,T^{*n}x\rangle=\Re \int_{\partial \omega}\overline{f(\sigma )}^{n+1}\langle N(s)Tx,T^{*n}x\rangle\,ds-\kappa ^2\|T^{*n}x\|^2.
\end{align*}
Therefore, with $\displaystyle r_1:=\Big(\int_{\partial \omega}\langle N(s)u(\sigma ),u(\sigma )\rangle\,ds\Big)^{1/2}$,
\begin{align}
\kappa m_n{-}m_{n+1}&=(\kappa ^2{-}\kappa )\|T^{*n}x\|^2-
\Re \int_{\partial \omega}\overline{f(\sigma )}^n\langle N(s)u(\sigma ),T^{*n}x\rangle\,ds\nonumber\\
&\geq (\kappa ^2{-}\kappa )\|T^{*n}x\|^2- \sqrt{\rho}\, \|T^{*n}x\|\,r_1\label{eq1b}\\
&\geq \big(\sqrt{\kappa ^2{-}\kappa}\, \|T^{*n}x\|-\frac{\sqrt{\rho} \,r_1}{2\sqrt{\kappa ^2{-}\kappa }}\big)^2-\frac{\rho \,r_1^2}{4(\kappa ^2{-}\kappa) }\\&\geq -\frac{\rho \,r_1^2}{4(\kappa ^2{-}\kappa) }.\nonumber
\end{align}
This implies
\begin{equation}
m_1\geq \kappa ^{-k}m_{n+1}-\frac{\rho \,r_1^2}{4(\kappa ^2{-}\kappa) }\frac{1{-}\kappa ^{-k}}{\kappa {-}1}.
\label{m1}
\end{equation}
Also, $m_{n+1}$ is bounded since, if $\| g \|_{\infty , \Omega} \leq 1$, then
\[
\| E(g) \| \leq \frac{1}{2 \pi} \int_{\partial \Omega} \| ( \sigma I - A )^{-1} \| | d \sigma | +
\int_{\partial \omega} | \lambda (s) |\,ds := C_E .
\]
This constant is finite and independent of $g$ since the spectrum of $A$ lies inside $\Omega$ so the
resolvent is bounded on $\partial \Omega$, and $\lambda$ is bounded and measurable on the finite-length boundary, and hence integrable.
Since $\| f \|_{\infty , \Omega} \leq 1 \Rightarrow \| f^{n+1} \|_{\infty , \Omega} \leq 1$, it follows that
$\| E( f^{n+1} ) \| \leq C_E$ for all $n \in \N$.  Since $E( f^{n+1} ) = \int_{\partial \omega} 
\overline{ f( \sigma )}^{n+1} N(s)\,ds - T^{*^{n+1}}$ and $\| \int_{\partial \omega} \overline{f( \sigma )}^{n+1} N(s)\,ds \|
\leq \rho$, it follows that $\| T^{n+1} \| = \| T^{*^{n+1}} \| \leq \rho + C_E$.  Consequently,
\[
| m_{n+1} | = | \Re \langle E( f^{n+1} ) T^{n+1} x,x \rangle | \leq  C_E ( \rho + C_E ) .
\]
Since $m_{n+1}$ is bounded, we get, by taking $k\to\infty$ in (\ref{m1}), $\displaystyle m_1\geq -\frac{\rho  \,r_1^2}{4 \kappa(\kappa{-}1)^2 }$.

Recall that 
\[
m_1= \Re \langle E( f ) T x,x \rangle =\Re\int_{\partial \omega}\overline{f(\sigma )}\langle N(s)Tx,x\rangle\,ds -\|Tx\|^2,
\]
therefore we have
\begin{align}
r_1^2&= 
\int_{\partial \omega}\langle N(s)\overline{f(\sigma )}Tx,\overline{f(\sigma )}Tx\rangle\,ds+\kappa ^2\int_{\partial \omega}\langle N(s)x,x\rangle\,ds\nonumber
-2\kappa \Re\int_{\partial \omega}\overline{f(\sigma )}\langle N(s)Tx,x\rangle\,ds\label{eq2}\\
& \le \ \rho \,\|Tx\|^2+\rho\, \kappa ^2\|x\|^2-2\,\kappa \,\|Tx\|^2-2\,\kappa \,m_1\\
& = \ 2\,\kappa ^2(\rho {-}\kappa )-2\,\kappa \,m_1\leq  2\,\kappa ^2(\rho {-}\kappa )+\frac{\rho \,r_1^2}{2(\kappa {-}1)^2}.
\nonumber
\end{align}
This leads to
\[
r_1^2\big(1-\frac{\rho }{2(\kappa {-}1)^2}\big)\leq 2\,\kappa ^2(\rho{-}\kappa ),
\]
which is impossible if $\kappa >\max\big(\rho ,1{+}\sqrt{\rho/2 }\big)$.

\vspace{.1in}
\noindent
{\em Remark}:  It was already noted in \cite[Remark (v)]{schw} that this implies the optimal 
spectral constant $2$ for the quantum annulus.  See also \cite{crzx2b,jutsi,pas,tsik}.  
In this paper, we use it to derive spectral constants for a variety of non-convex sets.\medskip

\noindent{\em Remark}: If $W(A)$ is contained in a convex domain $\Omega$, we can choose $\lambda(s) =0$, which gives $K=2$. Letting $\Omega$ shrink to $W(A)$ recovers the sharp bound $2$ for $W(A)$. In the general case, a natural choice is to take for $\lambda (s)$ the minimum of the spectrum of $\mu (\sigma (s),A)$.

\section{Estimates of $\lambda(s)$}
Fix a point $\sigma_0 := \sigma ( s_0 ) \in \partial \Omega$,
where the unit tangent $\sigma_0' := \sigma' ( s_0 )$ exists.
The half-plane $\Pi_0 := \{ z \in \C : \Im ( \sigma_0' ( \overline{\sigma_0} -
\bar{z})) \geq 0 \}$ has the same outward normal as $\Omega$ at $\sigma_0$.
See Figure \ref{fig:Pi0}.
Let $Sp( \cdot )$ denote the spectrum and $w( \cdot )$ the numerical radius and
let $\lambda_{min} (s) := \min Sp ( \mu ( \sigma (s) , A ) )$.
The following results are proved in \cite[Lemmas 5-8]{CG2019}:

\begin{figure}[ht]
\centerline{\epsfig{file=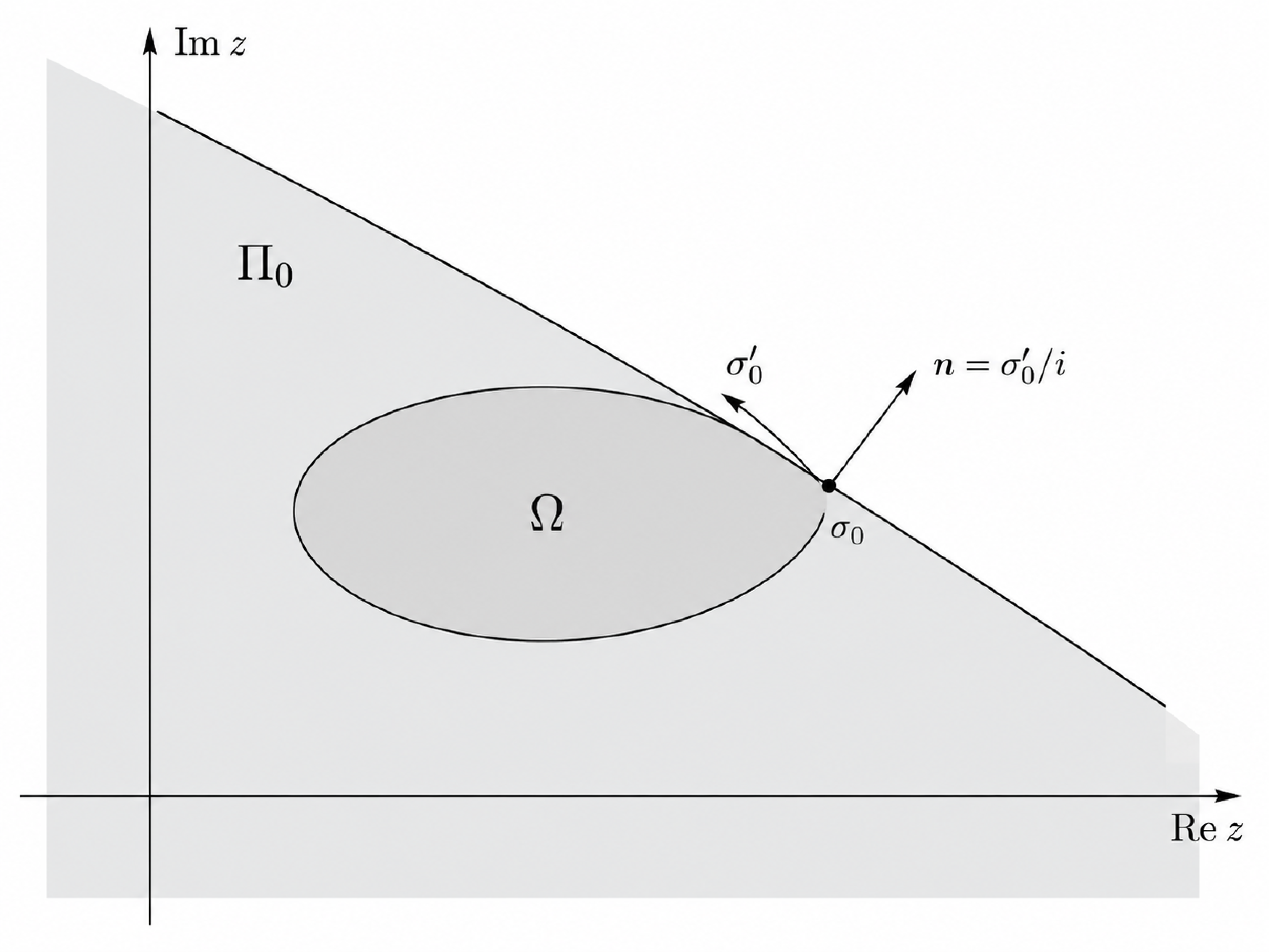,width=2in}}
\caption{$\Omega$ and the half-plane $\Pi_0$.}
\label{fig:Pi0} 
\end{figure}

\begin{lemma} \label{lem:lambdamin}

\begin{enumerate}
\item
If $W(A) \subset \Pi_0$, then $\lambda_{min} ( s_0 ) \geq 0$,
with equality if $\sigma_0 \in \partial W(A)$.
\item
Assume $| \sigma_0 - \xi | = R$ and $\{ z \in \C : | z - \xi | \leq R \} \subset \Pi_0$.
If $\| A - \xi I \| \leq R$, then $\lambda_{min} ( s_0 ) \geq \frac{1}{2 \pi R}$.
\item
If, for some $\xi \in \C \backslash Sp(A)$, $\sigma_0 - \xi = i r_1 \sigma_0'$,
where $0 < r_1 \leq 1 / \| (A- \xi I )^{-1} \|$, then $\lambda_{min} ( s_0 ) \geq - \frac{1}{2 \pi r_1}$.
\item
If, for some $\xi \in \C \backslash Sp(A)$,  $\sigma_0 - \xi = i r_2 \sigma_0'$,
where $0 < r_2 \leq 1/ w (( A - \xi I )^{-1} )$, then
$\lambda_{min} ( s_0 ) \geq - \frac{1}{ \pi r_2}$.
\end{enumerate}
\end{lemma}

\section{Examples}
\subsection{Annulus under norm assumptions}
\begin{figure}[ht]
\centerline{\epsfig{file=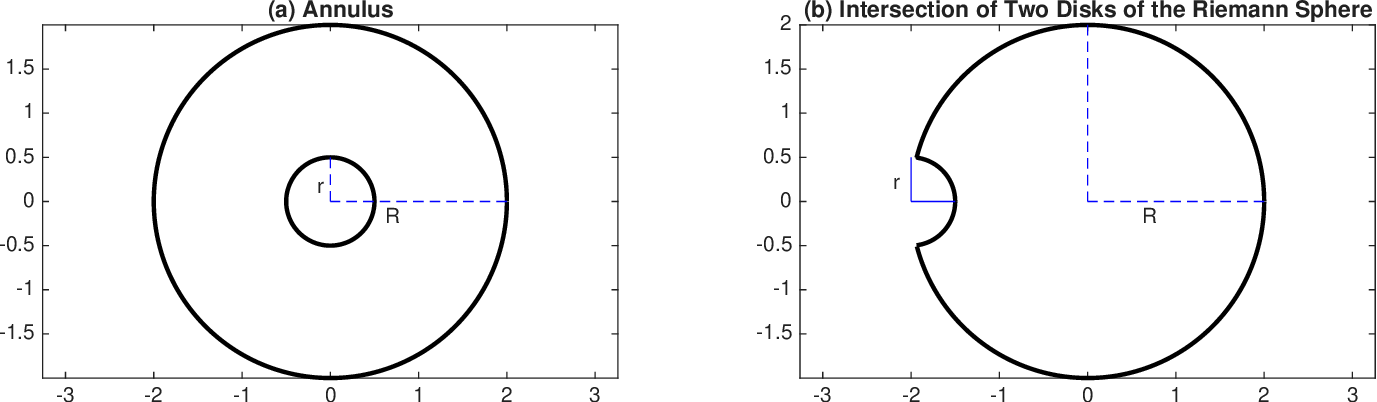,width=5in}}
\caption{(a) Annulus with outer radius $R > 1$ and inner radius $r = 1/R$.
(b) Intersection of disk of radius $R$ with exterior of disk of
radius $r$ centered at $(-R,0)$.}
\label{fig:annulus}
\end{figure}

Let ${\cal A}_R := \{ z \in \C : R^{-1} < |z| < R \}$, where $\| A \| < R$
and $\| A^{-1} \| < R$. Let $\Gamma_R = \{ z : |z| = R \}$ and $\Gamma_r =
\{ z : |z| = r \}$, where $r = R^{-1}$.  Assume that $\Gamma_R$ is oriented
counterclockwise and $\Gamma_r$ clockwise, in accordance with the positive orientation
of $\partial {\cal A}_R$.  See Figure \ref{fig:annulus} (a).

Define the modified double-layer kernel
\begin{equation}
N(s) = \mu ( \sigma (s) , A ) - \frac{1}{2 \pi i} \frac{\sigma' (s)}{\sigma (s)} I ,~~
\sigma \in \partial {\cal A}_R , \label{Pforannulus}
\end{equation}
where $\sigma'$ denotes differentiation with respect to arclength.
On $\Gamma_R$, $\sigma = R e^{i \theta}$ and $s = R \theta$, so
that $\sigma' = \frac{i \sigma}{R}$, and $N(s) = \mu( \sigma (s) , A) - \frac{1}{2 \pi R} I$.
Since $\| A \| < R$, Lemma \ref{lem:lambdamin}, part (2) shows that $N(s) \succeq 0$ for all
$\sigma \in \Gamma_{R}$.  On the inner circle, which has the opposite
orientation, $\sigma = r e^{-i \theta}$, $s = r \theta$,
so that $\sigma' = - \frac{i \sigma}{r} = -i R \sigma$ and
$N(s) = \mu ( \sigma (s), A) + \frac{R}{2 \pi} I$.
Since $r = 1/R < 1/ \| A^{-1} \|$, Lemma \ref{lem:lambdamin}, part(3)
shows that $N(s) \succeq 0$ for all $\sigma \in \Gamma_r$.
It follows that $\rho$ in Theorem \ref{th} satisfies
\[
\rho = 2 - \int_{\Gamma_R} \frac{1}{2 \pi R}\,ds + \int_{\Gamma_r} \frac{1}{2 \pi r}\,ds = 2 ,
\]
so that ${\cal A}_R$ is a $2$-spectral set for $A$.  This value $K=2$ is sharp \cite{tsik} and
 replaces the previous bound $K \leq 1 + \sqrt{2}$
in \cite{CG2019}.

{\bf Intersection of two disks of the Riemann sphere.}  See Figure \ref{fig:annulus} (b).
Theorem 10 of \cite{CG2019} reduces these regions by M\"{o}bius transformations either to
the annulus under norm assumptions case or to the convex domain containing $W(A)$ case.
That theorem can therefore be replaced by the following:

\begin{theorem}
Let $\Omega = D_1 \cap D_2$ with $D_1 = \{ z \in \C : | z - \omega_1 |
< R_1 \}$, $D_2 = \{ z \in \C : |z - \omega_2 | > 1/ R_2 \}$.
If $A$ satisfies $\| A - \omega_1 I \| \leq R_1$ and
$\| ( A - \omega_2 I )^{-1} \| \leq R_2$, then $\Omega$ is a
$2$-spectral set for $A$.
\end{theorem}

\subsection{Intersection of a bounded convex domain containing
$\overline{W(A)}$ with the exterior of a disk}

\begin{figure}[ht]
\centerline{\epsfig{file=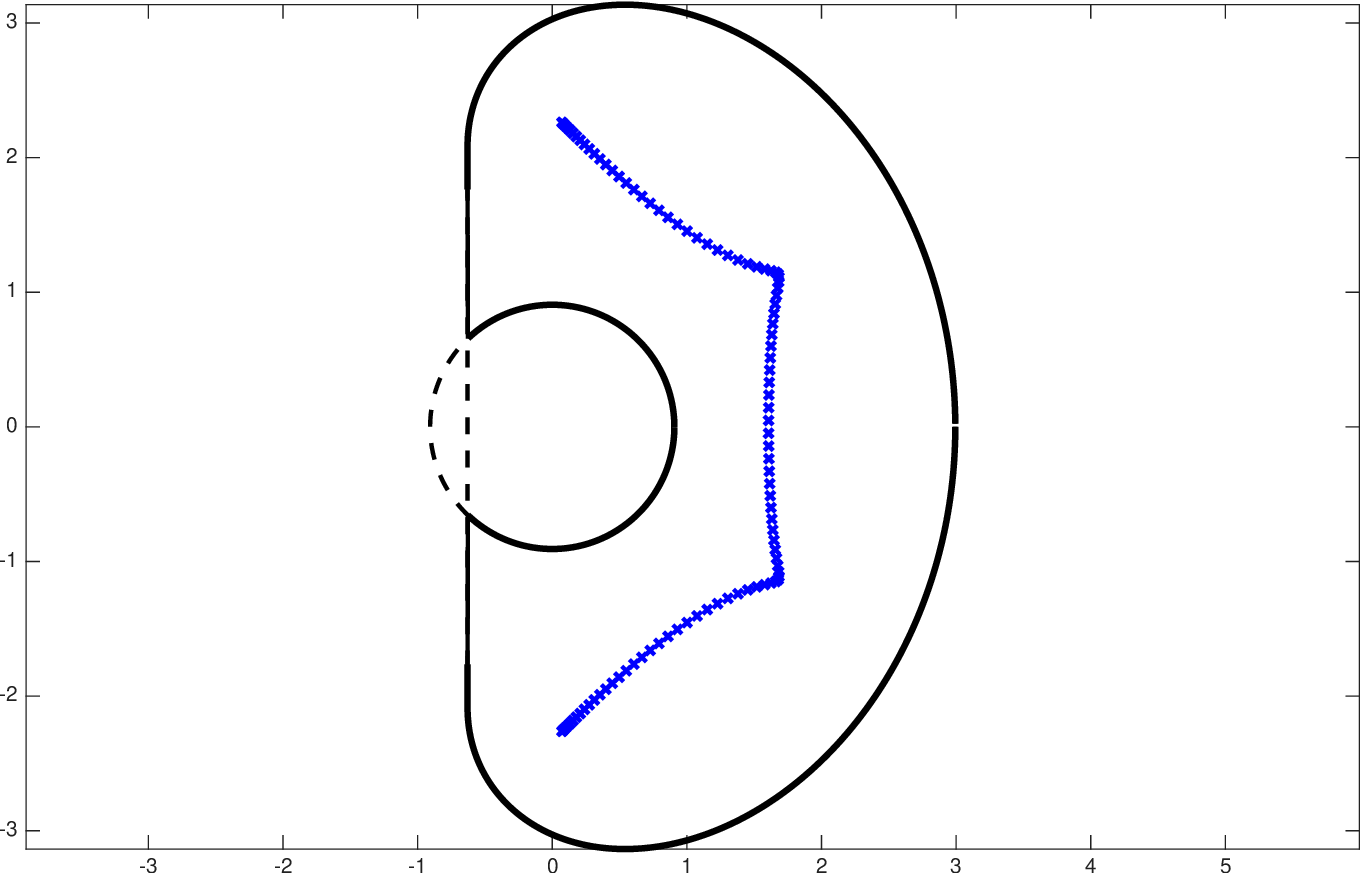,width=3in}}
\caption{Intersection of $W(A)$ with the exterior of a disk about
the origin.  Matrix is the Grcar matrix of order $100$.  Eigenvalues
are shown with x's.}
\label{fig:grcar100}
\end{figure}

Let $\Omega = \Omega_1 \cap \Omega_2$, where $\Omega_1$ is a bounded
convex domain containing $\overline{W(A)}$ and
$\Omega_2 = \{ z \in \C : | z - \xi | > 1/R \}$, where
$R > w( ( A - \xi I )^{-1} )$.  See Figure \ref{fig:grcar100}.

Write the boundary of $\Omega$ as $\Gamma_1 \cup \Gamma_2$, 
where $\Gamma_1 = \partial \Omega_1 \cap \overline{\Omega_2}$ and
$\Gamma_2 = \partial \Omega_2 \cap \overline{\Omega_1}$.  Then
$\mu ( \sigma (s) , A ) \succeq 0$ on $\Gamma_1$ and, by Lemma \ref{lem:lambdamin}, 
part (4),
$\nu ( \sigma (s),A ) := \mu ( \sigma (s) , A ) + \frac{R}{\pi} \succeq 0$
on $\Gamma_2$.
It follows that $\rho$ in Theorem \ref{th} satisfies
\[
\rho = 2 + \int_{\Gamma_2} \frac{R}{\pi}\,ds = 2 + \frac{R}{\pi} | \Gamma_2 | \leq 4 .
\]
Since $\rho \geq 2$, the first term in the maximum in Theorem \ref{th} dominates,
so that $\Omega$ is a $2 + \frac{R}{\pi} | \Gamma_2 | \leq 4$ spectral set for $A$.
This replaces the previous bound of $3 + 2 \sqrt{3}$ and retains
geometric information:  a smaller circular cutout gives a bound
closer to $2$.

\subsection{Annulus under numerical radius assumptions}
Applying the cutout construction to the full annulus under
the assumptions $w(A) < R$, $w( A^{-1} ) < R$, makes $\Gamma_2$
equal to the entire inner circle.
Hence $\frac{R}{\pi} | \Gamma_2 | = 2$ so $\rho = 4$ and $K \leq 4$.
This replaces the bounds in \cite[Theorem 11]{CG2019}, which range
between $6$ and $3 + \sqrt{10}$, depending on $R$.

\subsection{Removing multiple disks}
One can also improve on the $K$-values derived in \cite{GW2024} when multiple 
disks are removed from $\overline{W(A)}$.
Suppose $\Omega_0$ is a bounded convex domain containing $\overline{W(A)}$
and $\Omega$ is obtained
from $\Omega_0$ by removing $m$ disks centered at points
$\xi_1 , \ldots , \xi_m$, not in $\mbox{Sp}(A)$, with the radius $r_j$ of disk $j$
less than or equal to $1/w( (A - \xi_j I )^{-1} )$.  Set $p_j = 1$
if $r_j$ is also less than or equal to $1/ \| ( A - \xi_j I )^{-1} \|$;
otherwise, set $p_j = 2$.  Let $\Gamma_j =
\partial \Omega \cap \partial D( \xi_j , r_j )$.
On $\Gamma_j$, replace $\mu ( \sigma , A )$ by
\[
\mu ( \sigma , A ) + \frac{p_j}{2 \pi r_j} I ,
\]
which is positive semidefinite by Lemma \ref{lem:lambdamin}.
On the remaining part of $\partial \Omega_0$, $\mu$ is positive
semidefinite by convexity.  The resulting positive kernel
has total mass
\[
\rho I = \left( 2 + \sum_{j=1}^m \frac{p_j | \Gamma_j |}{2 \pi r_j} \right) I .
\]
Theorem \ref{th} then shows that $\Omega$ is a $K$-spectral set for $A$ with
\[
K \leq 2 + \sum_{j=1}^m \frac{p_j}{2 \pi r_j} | \Gamma_j | \leq
2 + \sum_{j=1}^m p_j .
\]

\subsection{Regions in which the numerical range of the resolvent is close to a disk about the origin}
In \cite{GKS2026}, the sets
\[
S_{\epsilon} (A) := \{ z \in \C : \tau_2 ( (A-zI )^{-1} ) / \tau_1 ( (A-zI )^{-1} ) < \epsilon \} ,~~
\epsilon \in (0,1) ,
\]
were studied.  Here $A$ is a square matrix and $\tau_1$ and $\tau_2$ are the largest and 
second largest singular values of the resolvent.  On the boundaries of these sets, where
$\tau_2 (( A- \sigma I )^{-1})/ \tau_1 ((A- \sigma I )^{-1}) = \epsilon$, if $\epsilon$ is very small,
then the resolvent is close to the rank one matrix $\tau_1 w_1 x_1^{*}$, where $w_1$ and $x_1$
are the left and right singular vectors corresponding to the largest singular value $\tau_1$
of the resolvent.  In this case, if also $| x_1^{*} w_1 |$ is small, then the numerical range 
of the resolvent is close to a disk about $\frac{1}{2} \tau_1 ( x_1^{*} w_1 )$
of radius $\tau_1 /2$ \cite{GW2024}, or, if  $| x_1^{*} w_1 |$ is sufficiently small, then the
numerical range of the resolvent is close to a disk about the origin of radius $\tau_1 / 2 =
\| ( \sigma I-A )^{-1} \| / 2$.
It was shown in \cite{GKS2026} that for highly nonnormal matrices and small values of $\epsilon$,
these two conditions -- $\tau_2 / \tau_1 < \epsilon$ and $| x_1^{*} w_1 | \approx 0$ -- often
hold over large regions, which, under proper scaling of $A$, look much like pseudospectra.
Now, the smallest eigenvalue of $\mu ( \sigma (s) , A )$ is $\frac{1}{\pi}$ times the smallest
imaginary part of points in $W( \sigma'(s) ( \sigma (s) I - A )^{-1} )$.  
Since $W( \sigma'(s) ( \sigma (s) I - A )^{-1})$ is approximately equal to a disk about the
origin of radius $\| ( \sigma (s) I - A )^{-1} \| /2$, it follows that
\[
\lambda_{min} ( s ) \approx - \frac{1}{2 \pi} \| ( \sigma (s) I-A )^{-1} \| .
\] 
In this case, $\rho$ in Theorem \ref{th} satisfies
\[
\rho = 2 - \int_{\partial S_{\epsilon} (A)} \lambda_{min} ( s ) \,ds \approx
2 + \frac{1}{2 \pi} \int_{\partial S_{\epsilon} (A)} \| ( \sigma (s) I - A )^{-1} \|\,ds .
\]
If $\rho > 2$, the first term in Theorem \ref{th} dominates and the new bound on $K$
in this limiting case is $2$ plus the elementary bound obtained by replacing the norm of the Cauchy
integral by the integral of the norm of the resolvent.  In this case, the new bound on $K$
is not helpful; a similar lack of improvement was observed in \cite{GW2024} for a bound on $K$ from \cite{CG2019}.

\section{Applications}
The updates to the applications section of \cite{CG2019} are
immediate.  For the GMRES algorithm for solving a linear system
$A {\bf x} = {\bf b}$, if ${\bf r_k} := {\bf b} - A {\bf x_k}$
is the residual at step $k$, then
\[
\frac{\| {\bf r_k} \|}{\| {\bf r_0} \|} \leq 2 \min \{ \sup_{z \in W(A)}
| p_k (z) | : p_k \in {\cal P}_k ,~p_k (0) = 1 \} ,
\]
where ${\cal P}_k$ is the set of polynomials of degree at most $k$.
This bound is not useful, however, if $W(A)$ contains the origin.

The cutout region described in section 4.2, consisting of the
intersection of $W(A)$ and the exterior of a disk about the
origin of radius $1/R$, where $R$ is the numerical radius of $A^{-1}$,
may provide a better bound.  If the origin lies inside $W(A)$
but its distance to the boundary of $W(A)$ is less than $1/R$,
this set excludes the origin and does not surround the origin,
so a polynomial with value $1$ at the origin may have magnitude
strictly less than $1$ throughout this region, which has now been
shown to be a $4$-spectral set for $A$.

In the rational Arnoldi algorithm for approximating $f(A) {\bf b}$,
the approximation at step $m$ is of the form $r_m (A) {\bf b}$,
where $r_m = p_{m-1} / q_{m-1}$ is a rational function with a
prescribed denominator $q_{m-1} \in {\cal P}_{m-1}$.  Letting
${\bf f_m^{RA}}$ denote the rational Arnoldi approximation at step $m$,
inequality (15) in \cite{CG2019} can be replaced by
\[
\| f(A) {\bf b} - {\bf f_m^{RA}} \| \leq 4 \| {\bf b} \| \min
\left\{ \sup_{z \in W(A)} | f(z) - r_m (z) | : r_m \in {\cal P}_{m-1}/
q_{m-1} \right\} .
\]
Other bounds on this difference given in \cite{CG2019}
are modified similarly to incorporate the new constants.

\vspace{.1in}
{\bf AI Disclosure Statement.} The second author used ChatGPT 5.6 Pro by OpenAI
in many aspects of this work.  The paper was written by the authors,
who take full responsibility for its content.


\begin{thebibliography}{11}

\bibitem{crzx2b} {\sc M.~Crouzeix},
{\em Spectral estimates in the quantum and in the numerical annulus}, 
 preprint, arXiv:2512.11813, 2025.

\bibitem{CG2019} M.~Crouzeix and A.~Greenbaum, {\em Spectral Sets:
Numerical range and beyond}, SIAM J.~Matrix Anal.~Appl., 40(3), 2019,
pp.~1087-1101.

\bibitem{GKS2026} A.~Greenbaum, F.~Kyanfar, and A.~Salemi,
{\em When is the resolvent like a rank one matrix?}, SIAM J.~Matrix Anal.~Appl.,
47(3), 2026, pp.~1326-1348.

\bibitem{GW2024} A.~Greenbaum and N.~Wellen,
{\em Comparison of $K$-spectral set bounds on norms of functions of a matrix
or operator}, Lin.~Alg.~Appl. 694 (2024), pp.~52-77.

\bibitem{jin} {\sc S.~Jin},
{\it The numerical range is a 2-spectral set}, 
Preprint on :https://www.preprints.org/manuscript/202607.1919, 2026.

\bibitem{jutsi} {\sc M.T.~Jury and G.~Tsikalas}, {\em Positivity conditions on the annulus via the double-layer potential kernel,} preprint, arXiv:2307.13387, 2023.

 \bibitem{schw} {\sc E.~Lorist and F.~Schwenninger},
{\it A solution to Crouzeix's conjecture}, arXiv:2608.03841v2, 2026.

 \bibitem{pas} {\sc J.E.~Pascoe}, {\em The spectral constant for the quantum cross and asymptotically sharp bounds for annuli},\ preprint, arXiv:2505.06230, 2025.
 
\bibitem{tsik}{\sc G.~Tsikalas}, {\em A note on a spectral constant associated with an annulus,} Operators \& Matrices, 16(1), (2022), pp.~95--99.

\end{thebibliography}
\end{document}